\documentclass[11pt,english,a4paper]{article}

\usepackage[utf8x]{inputenc}
\usepackage[T1]{fontenc}
\usepackage{cmlgc}
\usepackage{ucs}

\usepackage{xcolor}

\usepackage{graphicx}
\usepackage{mathtools,amsfonts,amssymb}
\usepackage{mathrsfs} 

\definecolor{vertfonce}{rgb}{0.20, 0.46, 0.25}
\definecolor{rougefonce}{rgb}{0.64, 0.09, 0.20}

\usepackage[breaklinks=true,
            colorlinks=true,
            linkcolor=rougefonce,
            citecolor=vertfonce]{hyperref}

\usepackage{autonum} 

\newcommand{\RM}{\mathbb{R}}
\newcommand{\ZM}{\mathbb{Z}}

\newcommand{\NM}{\mathbb{N}}
\newcommand{\CM}{\mathbb{C}}
\newcommand{\h}{\hbar}
\newcommand{\Cinf}{C^\infty}
\newcommand{\dd}[1]{\ensuremath{\operatorname{d}\!{#1}}}
\newcommand{\ham}[1]{\mathcal{X}_{#1}}
\newcommand{\theor}{Theorem}
\newtheorem{theo}{\theor}[section]
\newcommand{\defin}{Definition}
\newtheorem{defi}[theo]{\defin}

\newcommand{\restr}{\upharpoonright}

\newcommand{\Lh}{\mathcal{L}_\h}
\newcommand{\cI}{\mathcal{I}}

\title{Quantum footprints of Liouville integrable systems}

\author{\textsc{Vũ Ngọc} 
  San\footnote{Univ Rennes, CNRS, IRMAR - UMR 6625, F-35000 Rennes, France}}

\date{}

\begin{document}
\maketitle

\begin{abstract}
  We discuss the problem of recovering geometric objects from the
  spectrum of a quantum integrable system. In the case of one degree
  of freedom, precise results exist. In the general case, we report on
  the recent notion of good labellings of asymptotic lattices.
\end{abstract}

\begin{footnotesize}
  \noindent \textbf{MS Classification :}
  81S10. 
  81Q20, 
  58J50, 
  65L09. 
  \\
  To appear in \emph{Reviews in Mathematical Physics}, DOI: 10.1142/S0129055X20600144
\end{footnotesize}

\section{Inverse problems}

In this paper, we use for convenience the vocabulary of classical and
quantum mechanics, but one should keep in mind that inverse problems
can be stated in a more universal way. Our general question is: ``What
footprints does a classical system leave on its quantum counterparts,
and are they sufficient to recognize the classical system that
produced these footprints?''

Imagine you're walking in a snowy landscape, trying to take a good
photograph of a wild animal. A silhouette appears in the distance, you
have no idea what beast it can be; you seize your camera, look at the
small screen\dots, and the ghost has just disappeared. Just as if it
knew you wanted to capture it. The same goes for quantum particles,
they want to delocalize when they are observed. Yet we know they live
there. On the other hand, the footprints on the snow, they are real,
and stable. You can take your time and study them, until, maybe, by
clever induction, you find out what kind of animal was standing there.

What we have just described is an inverse problem: from the
observation of a signal emitted by some device, can we recognize the
device that has emitted the signal? If we hear the sound produced by
various instruments playing the same note C, can we tell the
instrument without looking? This question can easily be turned into a
mathematical problem, up to some simplifications, and was, for
instance, popularized by the ``Can one hear the shape of a drum''
paper by Kac~\cite{kac}. The ``sound'' is the superposition of all
possible ``frequencies'' of a drum, \emph{i.e.} the spectrum of the
Laplace operator on a Euclidean domain $\Omega\subset\RM^n$, and the
question is whether the shape of $\Omega$ (which is the equivalence
class of $\Omega$ under the action of the orthogonal group $O(n)$ and
translations) can be determined by the spectrum of the Dirichlet
Laplacian. Natural variants of this question exist. One can consider
Riemannian metrics on a compact manifold $X$: can one recover the
metric from the (discrete) spectrum of the corresponding
Laplace-Beltrami operator? Or, back to quantum mechanics, we may
consider a Schrödinger operator $-\h^2\Delta + V$ on $\RM^n$, and ask
whether its spectrum determines the electric potential $V$. Both
questions gave rise to an important literature, see for instance the
references cited in~\cite{zelditch-inverse-survey,
  datchev-hezari-survey13}
and~\cite{guillemin-wave,guillemin-paul-uribe}, and have applications
in non-quantum situations, for instance in seismology,
see~\cite{colin-inverse, dehoop-iantchenko-vdh-zhai-20}.

In this paper we consider the so-called ``symplectic case''. Given an
arbitrary ``semiclassical operator'' $\hat H_\h$ depending on a small
parameter $\h>0$ (see the next section) and its symbol $H$, which is a
smooth function on a symplectic manifold $(M,\omega)$, can you recover
the triple $(M,\omega,H)$ from the $\h$-family of spectra of
$\hat H_\h$, where $\h$ varies in a set accumulating at zero? This
natural question can be found for instance
in~\cite{iantchenko-sjostrand-zworski}.  Actually, Kac's problem, in
disguise, \emph{is} such a symplectic inverse problem; in this case
indeed, the semiclassical Laplacian is simply $\h^2\Delta$, and the
Hamiltonian $H$ is the metric on the cotangent bundle $M=T^* X$
induced by $g$.

\section{Quantization and Semiclassical analysis}

Semiclassical analysis is a general framework for obtaining a
``geometric limit'' from PDEs with highly oscillating solutions; the
name ``semiclassical'' comes from the model situation where classical
mechanics can be seen as a (singular) limit of quantum mechanics, as
Planck's constant $\h$ tends to zero. Our conventions in this text are
the following:

By \emph{Classical observables}, or classical Hamiltonians, we mean
smooth functions $H\in\Cinf(M)$ on a symplectic manifold $(M,\omega)$,
the classical phase space. For instance, $M= \RM^{2n}$ with the
canonical symplectic structure $\omega = \dd \xi \wedge \dd x$. Each
function $H$ defines an evolution equation, the flow of the associated
Hamiltonian vector field $\ham{H}$ defined by
$\iota_{\ham{H}}\omega = -\dd H$.

By \emph{Quantum observables}, or quantum Hamiltonians, we mean a
selfadjoint operator $\hat{H}$ on a Hilbert space $\mathcal{H}$, and
the Hilbert space itself must be the quantization of a classical phase
space $M$. Each quantum Hamiltonian gives rise to the evolution
governed by the Schrödinger equation
\[
  \frac{\h}i\partial_t \psi = \hat H \psi,
\]
so, formally, $\psi(t) = e^{\frac{it}\h\hat H} \psi(0)$.  Stationary
states are solutions of the form
$\psi(t) = e^{\frac{i\lambda t}\h} u$, where $u\in\mathcal{H}$ and
$\hat{H}u = \lambda u$. It is a fascinating subject to understand the
relationships between the classical Hamiltonian flow and the quantum
Schrödinger evolution.

Two rigorous quantization schemes allow us to realise the above
picture: when $M$ is a cotangent bundle, $M=T^*X$, one can use
pseudo-differential quantization, see for
instance~\cite{martinez-book}. When $M$ is a prequantizable Kähler
manifold, one can use Berezin-Toeplitz quantization, see for
instance~\cite{lefloch-book} (Berezin-Toeplitz quantization was later
extended to general symplectic manifolds,
see~\cite{charles16}). Berezin-Toeplitz and pseudo-differential
quantization are similar in many respects, and both benefit from the
power and flexibility of microlocal analysis \emph{à la} Maslov,
Hörmander, etc..

\section{1D Hamiltonians : the Morse case}

Let $(M,\omega)$ be a 2-dimensional symplectic manifold.  Let
$H: M\to \RM$ a proper Morse function. Following the usual Morse
approach we will be interested in the (singular) foliation of $M$ by
level sets of $H$. An important object, the \emph{Reeb graph}
$\mathcal{G}$, is the set of leaves, \emph{i.e.} connected components
of level sets of $H$; in a neighbourhood of a regular level set,
$\mathcal{G}$ is a smooth, one-dimensional manifold. The smooth parts
are the edges of the graph. Each critical point of $H$ contributes to
a graph vertex, its degree is one for elliptic singularities (the
vertex is then a leaf), and three for hyperbolic singularities.

The Reeb graph turns out to be essential in the description of the
spectrum of a 1D quantum Hamiltonian. Let
$\hat H:=\textup{Op}_\h(H_\h)$ be the quantization of a symbol
$H_h:=H+\h H_1+\h^2H_2+\cdots$ on $M$.  Let $I\subset\RM$ be a closed,
bounded interval. Since $H$ is proper, $H^{-1}(I)$ is compact. In the
case of pseudo-differential quantization, we assume that $H_\h$
belongs to a symbol class and is elliptic at infinity, see for
instance~\cite{martinez-book}. Then the spectrum of $\hat H$ in any
compact subset of $\mathop{\textup{int}} I$ is discrete. In this case,
the inverse spectral theory is well understood, and summarized by the
following statement.

\begin{theo}[\cite{san-inverse,leflo2014}]
  Let $M_I := H^{-1}(\mathop{\textup{int}} I)$.  Suppose that
  $H_{\restr {M_I}}$ is a simple Morse function.  Assume that the
  graphs of the periods of all trajectories of the Hamiltonian flow
  defined by $H_{\restr {M_I}} $, as functions of the energy,
  intersect generically.

  Then the knowledge of the spectrum
  $\sigma(\hat H)\cap I+\mathcal{O}(\h^2)$ determines
  $(M_I,\omega,H)$.
\end{theo}

The proof of this theorem, like many of its kind, is divided in two
steps.  The first one is to recover the Reeb graph $\mathcal{G}$ of
$(M_{I},H_I)$ from the spectrum. The second step is to prove that
$\mathcal{G}$, decorated with appropriate numerical invariants that we
can also recover from the spectrum, completely determines the
classical system $(M_I,\omega,H)$.  The last step was proven by
Dufour-Toulet-Molino~\cite{dufour-mol-toul}. The first step was
established in the pseudo-differential case in ~\cite{san-inverse},
and in the case of Berezin-Toeplitz operators by Le
Floch~\cite{leflo2014}. It involved microlocal analysis in the
(time/energy) phase space to be able to separate the various connected
components of $\mathcal{G}$ contributing to the same region of the
spectrum.

More recently, a new interpretation of this result has been proposed
by several mathematicians, in particular Leonid Polterovich and the
author. Suppose you add a generic \emph{non-selfadjoint} perturbation
to the quantum operator $\hat H$. Then, the connected components of
$\mathcal{G}$, instead of leading to overlapping parts of the spectrum
--- and hence potentially difficult to tell apart --- should instead
give rise to different \emph{complex branches} of the spectrum of the
non-selfadjoint operator. Thanks to the recent result by
Rouby~\cite{rouby-17} explaining the non-selfadjoint version of
Bohr-Sommerfeld quantization conditions, we believe that this
conjectural interpretation should produce new rigorous results.

\begin{theo}[Rouby~\cite{rouby-17}]
  Let $P_\epsilon$ be an \emph{analytic} pseudodifferential operator
  on $\RM$ or $S^1$ of the form $P_\epsilon = \hat H + i\epsilon Q$,
  where $\hat H$ is selfadjoint with discrete spectrum, and $Q$ is
  $\hat H$-bounded.

  Then, near any regular value of the symbol $H$, with connected
  fibers, the spectrum of $P_\epsilon$ is given by
  $\{g(\h m;\epsilon); m\in\ZM\}$, where
  $g:\CM\to\CM$ is holomorphic and\\
  \[
    g \sim g_0 + \h g_1 + \h^2 g_2 + \cdots
  \]
  Moreover, $g_0$ is the inverse of the action variable, and\\
  \[
    g_0 \sim H + i\epsilon\langle q \rangle + \mathcal{O}(\epsilon^2).
  \]
\end{theo}

Rouby's theorem is technically quite involved, because one has to take
advantage of analyticity to fight non-selfadjoint instability
(pseudo-spectral effects), and usual $\Cinf$ microlocal analysis is
not strong enough for this. No analogue of this result for
Berezin-Toeplitz quantization exists yet. However, very recent
advances on the analyticity of the Bergman projection give some hope,
see~\cite{san-rouby-sjostrand18, deleporte18, charles2019analytic}.
  
\section{Integrable systems}

In view of Rouby's theorem, one can notice that a particular case
where analyticity is not required occurs when $P_\epsilon$ is
\emph{normal}, \emph{i.e.} the non-selfadjoint perturbation $Q$
\emph{commutes} with the selfadjoint part $\hat H$. More generally, a
number of results exist in the presence of a \emph{completely
  integrable quantum system}, by which we mean the data of $n$
pairwise commuting selfadjoint operators $P_1,\dots, P_n$, when the
phase space $M$ is $2n$-dimensional. In fact, even for operators that
are not quantum integrable but still have a completely integrable
classical limit, quite precise results can be obtained, for both
direct and inverse problems; see~\cite{san-hitrik-sjostrand, hall-13},
and references therein.

This notion of quantum integrability parallels the usual
\emph{Liouville integrability} of classical Hamiltonians, where we
dispose of $n$ independent Poisson-commuting functions $f_1,\dots,f_n$
on $M$. Note that, near a \emph{regular} level set of the joint map
$F:=(f_1,\dots,f_n): M \to \RM^n$, one has \emph{action-angle}
coordinates, due to the celebrated Liouville-Mineur-Arnold theorem,
but Liouville integrability is \emph{more general}: it allows for
singularities where the action-angle theorem cannot apply.

The natural multi-dimensional generalization of the Reeb graph is the
leaf space of the ``moment map'' $F$, which is equipped with a natural
\emph{integral affine structure} (see for
instance~\cite{san-sepe-poisson16}). The quantum analogue of this
singular integral affine manifold is the \emph{joint spectrum} of the
commuting operators $P_1,\dots, P_n$. Hence, we are naturally lead to
the following inverse problem: given an $\h$-family of joint spectra,
can one recover the triple $(M,\omega,F)$?

A first approach to this question is to restrict oneself to
Hamiltonian systems with many compact symmetries; namely the
\emph{toric} and \emph{semitoric}
cases. See~\cite{san-alvaro-first-steps} for a description of a
general program of study, and conjectures.

\section{Asymptotic lattices}

Having in mind the general inverse problem for quantum integrable
systems, another angle of attack is to consider the regular part of
the integral affine structure, and exploit the lattice structure of
the joint spectrum, which was already established by Colin de
Verdière~\cite{colinII}. This leads to the notion of \emph{asymptotic
  lattices}, whose systematic study was initiated recently
in~\cite{san-dauge-hall-rotation}. Although the initial goal of that
paper was to recover from the quantum spectrum a specific classical
invariant, \emph{the rotation number}, we believe that the general
setup should help understanding all invariants related to the integral
affine structure. In particular, we hope that it will allow to finally
obtain a complete result on the inverse theory of semitoric systems.

Let $B\subset\RM^n$ be a simply connected bounded open set. Let
$\Lh\subset B$ be a discrete subset of $B$ depending on the small
parameter $\h\in\cI$, where $\cI\subset\RM^*_+$ is a set of positive
real numbers admitting $0$ as an accumulation point. Here is a
slightly imprecise definition of asymptotic lattices (we don't delve
into multiplicity issues and the details of the $\O(\h^\infty)$
topology).

\begin{defi}[Asymptotic lattice~\cite{san-mono,
    san-dauge-hall-rotation}]
  We say that $(\Lh, \cI, B)$ is an \emph{asymptotic lattice} if
  \[
    \Lh = G_\h(\h\ZM^n \cap U) + \mathcal{O}(\h^\infty)
  \]
  with
  \[
    G_\h = G_0 + \h G_1 + \h^2 G_2 + \cdots
  \]
  in the $\Cinf(U)$ topology, where $G_0:U\to\RM^n$ is a
  diffeomorphism on its image.
\end{defi}

The definition is motivated by the following older result.
\begin{theo}[\cite{colinII,charbonnel,charles-bs}]
  Let $P:=(P_1,\dots, P_n)$ be a Quantum integrable system. Let
  $c\in\RM^n$ be a regular value of the classical moment map $F$ with
  connected fiber. Then the \emph{joint spectrum} of $P$ near $c$ is
  an \emph{asymptotic lattice}.
\end{theo}

\section{Good labellings}

In order to recover the integral affine structure from the spectrum,
one needs to recover the map $G_0$ in the previous theorem. In order
to do this, we claim that it is enough to find a ``good labelling'' of
all joint eigenvalues. By this we mean, to assign to each joint
eigenvalue $\lambda$ a $n$-uple of integers $(k_1,\dots,k_n)$ such
that
\[
  \lambda = G_\h(\h k_1,\dots, \h k_n) + \mathcal{O}(\h^\infty).
\]

In~\cite{san-dauge-hall-rotation}, we investigated the case of two
degrees of freedom, $n=2$.  To our surprise, this question turned out
to be more intricate (and more interesting) than what we initially
thought. 

On the other hand, the process of finding a good labelling is
elementary, and can be described algorithmically, which is important
for the following reason. The informal question ``can one hear the
shape of a drum'' has two possible interpretations. The minimal one is
to prove \emph{injectivity} of the map sending a classical system to
its quantum spectrum. In this case, the classical system is
\emph{determined} by the quantum spectrum in a weak sense: two
different classical systems cannot give rise to the same quantum
spectrum. A stronger result would be to obtain the classical system
that produced the spectrum in a constructive way. Writing an algorithm
contributes to the latter.

The algorithm will be constructed in two steps. In the first one, the
value $\h$ is fixed, and the algorithm returns a candidate labelling
$\lambda\mapsto (k_1,k_2)$. However, this candidate does not have the
required continuity property in the variable $\h$. Hence we perform a
second step where we consider now a full sequence of values
$\h_i\to 0$, and we correct the discontinuity by an inductive
algorithm in the variable $i\in\NM$. We don't know whether a direct
approach, in one step, would be possible. When constructing a good
labelling, another difficulty comes from the choice of a valid
``origin'' for the lattice. For this purpose, the set of values of
$\h$ must be ``dense enough'' when accumulating at zero. Values of the
form $\h=\frac{1}{k}$ for $k\in\NM\setminus\{0\}$ do not fulfill this
requirement, which is an issue for Berezin-Toeplitz
quantization. However, in many applications, the choice of the lattice
origin is irrelevant. Introducing the notion of ``linear labelling''
as a good labelling ``modulo its origin'', we have the following
result.

\begin{theo}
  \label{theo:algo} There exists an explicit algorithm such that the
  following holds. Let $(\Lh, \cI, B)$ be an asymptotic lattice, where
  $B\subset\RM^2$. Let $\h_j\in\cI$, $j\geq1$, be a decreasing
  sequence tending to $0$.  Then, from this data, the algorithm
  produces a \emph{linear labelling
  } of the asymptotic lattice
  $(\Lh, \cI', B)$, where $\cI'=\{\h_j, j\in\NM^*\}$.
\end{theo}

\begin{figure}
  \centering \includegraphics[width=0.3\linewidth]{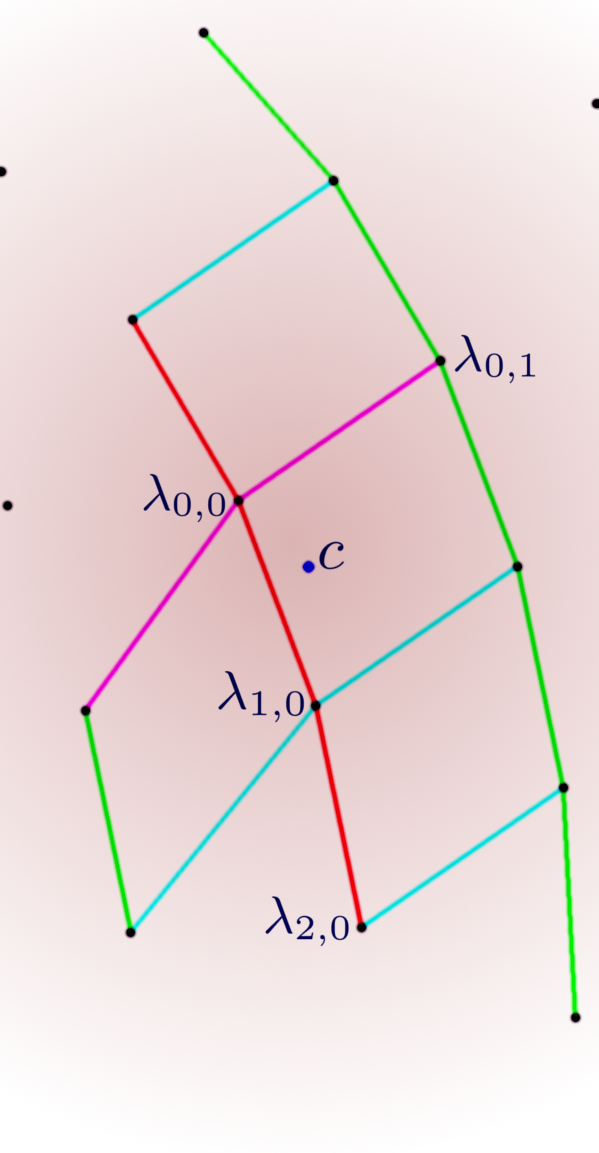}
  \caption{The labelling algorithm}
  \label{fig:algo}
\end{figure}

Below is the complete description of the first step. If $(k_1,k_2)$ is
a label for a point $\lambda\in\Lh$, we shall denote this point
$\lambda=\lambda_{k_1,k_2}$. The complete algorithm is as
follows~\cite{san-dauge-hall-rotation}, and pictured in
Figure~\ref{fig:algo}.
\begin{enumerate}
\item \label{item:algo-c} Choose an open subset $B_0\Subset B$, and
  fix $c\in B_0$.
\item Choose a closest point to $c$. Label it as $(0,0)$.

\item \label{item:algo-1-0} Choose a closest point to $\lambda_{0,0}$
  (in the set $\Lh\setminus \{\lambda_{0,0}\}$).  Label it as $(1,0)$.

\item \label{item:algo-n-0}

  Choose a closest point to
  $2\lambda_{1,0} - \lambda_{0,0} =\lambda_{1,0} +
  (\lambda_{1,0}-\lambda_{0,0})$ and label it as $(2,0)$.

  Continuing in this fashion (if $\lambda_{k_1-1,0}$ is chosen, take
  $\lambda_{k_1,0}$ to be a closest point to
  $\lambda_{k_1-1,0} +(\lambda_{k_1-1,0} - \lambda_{k_1-2,0})$), label
  points $\lambda_{k_1,0}$, $k_1>0$, until the next point lies outside
  of $B_0$.

  Label $\lambda_{k_1,0}$ for negative $k_1$ in the same way, starting
  from the closest point to $2\lambda_{0,0} - \lambda_{1,0}$.

\item \label{item:algo-0-1} Choose a closest point to $\lambda_{0,0}$
  not already labeled and label it as $(0,1)$.

\item \label{item:algo-1-1} Label a closest point to
  $\lambda_{0,1} + (\lambda_{1,0} - \lambda_{0,0})$ as $(1,1)$.

\item Use the points $\lambda_{0,1}$, $\lambda_{1,1}$ to repeat the
  process in step 4, labelling as many points $\lambda_{k_1,1}$ as
  possible (if $\lambda_{k_1-1,1}$ is chosen, take $\lambda_{k_1,1}$
  to be a closest point to
  $\lambda_{k_1-1,1} +(\lambda_{k_1-1,1} - \lambda_{k_1-2,1})$).

\item \label{item:algo-0-2} Label a closest point to
  $2\lambda_{0,1} - \lambda_{0,0}$ as $\lambda_{0,2}$. Repeat steps
  6-7 to label all points $\lambda_{k_1,2}$.

\item Continuing as above, label all points $\lambda_{k_1,k_2}$,
  $k_2>0$ which lie in the given neighborhood.

\item Label a closest point to $2\lambda_{0,0} - \lambda_{0,1}$ as
  $(0,-1)$.
\item Repeat steps 6,7,8,9 with negative $k_2$ indices.

\item\label{item:algo-determinant} Finally, if the determinant of the
  vectors
  $(\lambda_{1,0} - \lambda_{0,0},\lambda_{0,1} - \lambda_{0,0})$ is
  negative, switch the labelling
  $\lambda_{k_1,k_2} \mapsto \lambda_{-k_1,k_2}$ (in order to make it
  oriented).
\end{enumerate}

This algorithm should be useful in several inverse spectral
problems. For instance, in~\cite{san-dauge-hall-rotation},
Theorem~\ref{theo:algo} was used to prove that the classical
\emph{rotation number} of any Liouville torus can be recovered from
the joint spectrum. From a quite different perspective, it could also
be interesting to investigate the proximity of our approach with
\emph{topological data analysis} and \emph{manifold learning}.

\section{Prospects}

The detection of good labellings should allow the complete
reconstruction of the integral affine structure, at least for its
regular part.  The next step would be to \emph{globalize} the notion
of asymptotic lattice, and include \emph{singularities}, to obtain
\emph{quantized integral affine structures with singularities}. For
instance, the singular limit of the rotation number, as explained
in~\cite{san-dullin}, should be a feature of asymptotic lattices with
focus-focus singularities.

As far as inverse spectral theory is concerned, we certainly hope to
use the notion of asymptotic lattice to advance towards the
\emph{Spectral semitoric conjecture}
\cite{san-alvaro-survey,san-lefloch-pelayo:jc}: can you detect the
five symplectic invariants of a semitoric system on the joint spectrum
(or asymptotic lattice)? In an ongoing work with Le Floch, which
intially focussed on the reconstruction of the \emph{twisting index}
invariant, we finally expect to obtain not only the injectivity of the
``semiclassical joint spectrum map'' for simple semitoric systems, but
also a full reconstruction procedure.  A more general result should
include multi-pinched tori
(see~\cite{pelayo-tang-conjecture,palmer-pelayo-tang2019semitoric,
  meulenaere-hohloch2019}), and systems with non-proper circle moment
map~\cite{san-pelayo-ratiu-connectivity,san-pelayo-ratiu-affine}.

\bibliographystyle{abbrv}%
\bibliography{bibli-utf8}
\end{document}